\documentclass[10pt, twoside, a4paper]{article}
\let\finishall\relax\let\Finishall\relax\let\getprepared\relax
\ifx\TestIngCommand\undefined\relax
\else\input ../../../proc/ppreamb-book.tex 
  \input init.tex \fi
\let\TestIngCommand\undefined

\newtheorem{remark}{Remark}

\newtheorem{conj}{Conjecture}
\usepackage{amsmath,amscd,amssymb}                                         
\usepackage{graphics}                                                     
\newtheorem{theo}{Theorem}                                                 
\newtheorem{lem}{Lemma}                                                    
\newtheorem{defi}{Definition}
\newskip\ttglue\ttglue=.5em plus.25em minus.15em                           
\def\firstname#1{\def\FIRSTNAME{#1}\ignorespaces}
\def\lastname#1{\def\LASTNAME{#1}\ignorespaces}
\def\middleinitial#1{\def\MIDDLEINI{#1}\ignorespaces}
\def\department#1{\def\DEPARTMENT{#1}\ignorespaces}
\def\institute#1{\def\INSTITUTE{#1}\ignorespaces}
\def\address#1{\def\ADDRESS{#1}\ignorespaces}
\def\country#1{\def\COUNTRY{#1}\ignorespaces}
\def\otheraffiliation#1{\def\OTHERAFFILIATION{#1}\ignorespaces}
\def\email#1{\def\EMAIL{#1}\ignorespaces}
\newcount\autcount\autcount=0                                              
\newcount\affcount\affcount=0                                              
\newcount\numcount\newcount\nummcount                                      
\newcount\nummmcount\newcount\nummmmcount                                  
\def\writename#1#2{\ \kern-1ex\hbox{
  \csname AUthor\the#1\endcsname\                                          
  \edef\TESTSTR{}\expandafter\ifx\csname auTHor\the#1\endcsname\TESTSTR    
  \else\csname auTHor\the#1\endcsname.\ \fi                                
  \csname authOR\the#1\endcsname$^{\csname AFF\the#1\endcsname}$
  \expandafter\ifx\csname corr\number#1\endcsname\relax                    
  \else\thanks{Corresponding author.}\ \fi                                 
  }\ifnum#1<#2, \else\ \kern-1ex\fi}
\def\writeemail#1{
  \nummcount=0\relax\nummmcount=0\relax                                    
  \loop\ifnum\nummcount<\autcount\advance\nummcount by1\relax              
    {\expandafter\ifnum\csname AFF\the\nummcount\endcsname=#1\relax        
    \global\advance\nummmcount by1\fi}\repeat                              
  \nummcount=0\relax\nummmmcount=0\relax                                   
  \loop\ifnum\nummcount<\autcount\advance\nummcount by1\relax              
    {\expandafter\ifnum\csname AFF\the\nummcount\endcsname=#1\relax        
    \global\advance\nummmmcount by1\relax\def\blank{}\expandafter          
    \ifx\csname EMAIL\the\nummcount\endcsname\blank(no e-mail)
    \else\csname EMAIL\the\nummcount\endcsname                             
    \fi                                                                    
    \ifnum\nummmmcount<\nummmcount; \fi\fi}\repeat}
\long\def\BeginAuthorList#1\EndAuthorList{#1\relax                         
  \author{\vbox{\hsize=390pt\noindent\numcount=0\relax                     
    \loop\ifnum\numcount<\autcount\advance\numcount by1\relax              
      \writename{\numcount}{\autcount}
      \repeat}\\[2mm]                                                      
    \vbox{\small\numcount=0\relax                                          
      \loop\ifnum\numcount<\affcount\advance\numcount by1\relax            
        \vbox{{\count0=\numcount\relax                                     
          \loop\expandafter\ifnum\csname AFF\the\count0\endcsname
            <\numcount\relax\advance\count0 by1\relax\repeat               
          $^{\csname AFF\the\count0\endcsname}$}
        \def\BLANK{}\expandafter\ifx\csname DEPT\the\numcount\endcsname    
          \BLANK                                                           
          \else\csname DEPT\the\numcount\endcsname, \fi                    
        \csname INST\the\numcount\endcsname,                               
        \csname ADDR\the\numcount\endcsname,                               
        \csname COUN\the\numcount\endcsname                                
        \edef\TEST{}\expandafter\ifx\csname OTHE\the\numcount\endcsname
          \TEST                                                            
          .\else;\break\csname OTHE\the\numcount\endcsname.\fi}
        \vbox{\writeemail{\numcount}}
        \repeat}\\}}
\expandafter\def\csname x1\endcsname{}
\expandafter\def\csname x2\endcsname{}
\expandafter\def\csname x3\endcsname{}
\expandafter\def\csname x4\endcsname{}
\expandafter\def\csname x5\endcsname{}
\expandafter\def\csname x6\endcsname{}
\expandafter\def\csname x7\endcsname{}
\expandafter\def\csname x8\endcsname{}
\expandafter\def\csname x9\endcsname{}
\def\Author#1#2{\global\advance\autcount by1\relax#2                       
  \expandafter\edef\csname AUthor\the\autcount\endcsname{\FIRSTNAME}
  \expandafter\edef\csname auTHor\the\autcount\endcsname{\MIDDLEINI}
  \expandafter\edef\csname authOR\the\autcount\endcsname{\LASTNAME}
  \expandafter\edef\csname EMAIL\the\autcount\endcsname{\EMAIL}
  \let\tempera\"\def\"{\string\"}\expandafter\ifx\csname x\DEPARTMENT
    \endcsname\relax                                                       
    \global\advance\affcount by1\relax\let\"\tempera                       
    \expandafter\edef\csname DEPT\the\affcount\endcsname{\DEPARTMENT}
    \expandafter\edef\csname INST\the\affcount\endcsname{\INSTITUTE}
    \expandafter\edef\csname ADDR\the\affcount\endcsname{\ADDRESS}
    \expandafter\edef\csname COUN\the\affcount\endcsname{\COUNTRY}
    \expandafter\edef\csname OTHE\the\affcount\endcsname{\OTHERAFFILIATION}
    \expandafter\edef\csname AFF\the\autcount\endcsname{\the\affcount}
  \else\expandafter\edef\csname AFF\the\autcount\endcsname{\DEPARTMENT}
  \fi\let\"\tempera\ignorespaces}
\def\CorrespondingAuthor#1#2{
  \expandafter\xdef\csname corr\number#1\endcsname{cor}
  \Author#1{#2}}
\def\PaperTitle#1{\title{\bf#1}}
\def\Category#1{\ignorespaces}
\def\keywords#1{{\noindent \emph{Keywords:}                                
  \def\BLANK{}\def\TEST{#1}\ifx\BLANK\TEST(n/a).\else#1\fi}}
\getprepared                                                               
\begin{document}                                                           

\PaperTitle{Multiresolution Analyses on Quasilattices}
\Category{(Pure) Mathematics}

\date{}

\BeginAuthorList
  \Author1{
    \firstname{Wayne}
    \lastname{Lawton}
    \middleinitial{M}   
    \department{Applied Mathematics Program}
    \institute{Mahidol University}
    \address{Nakornpathom}
    \country{Thailand}
    \otheraffiliation{School of Mathematics and Statistics, University of Western Australia, Perth, Australia}
    \email{wayne.lawton@uwa.edu.au}}
\EndAuthorList
\maketitle
\thispagestyle{empty}
\begin{abstract}
We derive relations between geometric means of the Fourier moduli of a refinable
distribution and of a related polynomial. We use Pisot-Vijayaraghavan numbers to construct families of one dimension quasilattices and multiresolution analyses spanned by distributions that are refinable with respect to dilation by the PV numbers and translation by quasilattice points. We conjecture that scalar valued refinable distributions are never integrable, construct piecewise constant vector valued refinable functions, and discuss multidimensional extensions.
\end{abstract}

\finishall
\section{Refinable Distributions}
\label{intro}
In this paper $\mathbb{Z}, \mathbb{N} = \{1,2,3,...\}, \mathbb{Q}, \mathbb{R}, \mathbb{C}$ denote the integer, natural, rational, real, and complex numbers and $\mathbb{T} = \mathbb{R}/\mathbb{Z}$ denotes the circle group.
\\ \\
Let $\lambda \in \mathbb{R}\, \backslash \, [-1,1],$ $a_1,....a_m \in \mathbb{C}\backslash\{0\},$
$a_1+\cdots+a_m = |\lambda|,$ $\tau_1 < \cdots < \tau_m \in \mathbb{R},$ and
\begin{equation}
\label{A}
A(y) = |\lambda|^{-1} \sum_{j=1}^{m} a_j\, e^{2\pi i \tau_j y}, \ \ y \in \mathbb{R}.
\end{equation}
There is a unique compactly supported distribution $f$
satisfying the refinement equation
\begin{equation}
\label{refinable}
    f(x) = \sum_{j=1}^m a_j f(\lambda \, x - \tau_j)
\end{equation}
and $\int f(x)dx = 1.$ Its Fourier transform ${\widehat f}(y) = \int_{-\infty}^{\infty} \, f(x)e^{2\pi i x y} \, dx$ has the expansion
\begin{equation}
\label{product}
    {\widehat f}\,(y) = \prod_{k=1}^{\infty} A\left(\, y\lambda^{-k}\, \right)
\end{equation}
and ${\widehat f}\,(0) = 1.$ Therefore ${\widehat f}$ and the polynomial $A$ are related by the functional equation
\begin{equation}
\label{functional}
    {\widehat f}\, \left(\, y \, \lambda^k \, \right) = {\widehat f}\,(y) \,
    \prod_{j = 0}^{k-1} A\left(\, y \, \lambda^j\, \right), \ \ k \in \mathbb{N}.
\end{equation}
The Mahler measure \cite{mahler} (or height) $M(P)$ of a polynomial $P(z_1,...,z_d)$ is defined by
\begin{equation}
\label{mahler}
    M(P) = \exp \int_{0}^{1} \cdots \int_{0}^{1}
            \ln |P(e^{2\pi i \, t_1},...,e^{2\pi i \, t_d})|\, dt_1 \dots dt_d.
\end{equation}
For $P(z) = c(z-\omega_1)\cdots(z-\omega_q)$ Jensen's theorem \cite{jensen} gives
$M(P) = |c| \prod_{j = 1}^{q} \max\,  \{\, 1,|\,\omega_j\,|\, \}.$
\\ \\
\begin{remark} Mahler measure arises in prediction theory of stationary
random processes \cite{lawton1}, algebraic dynamics \cite{everestward}, and
in Lehmer's open problem in number theory \cite{lehmer}.
\end{remark}
$A$ is almost periodic in the (uniform) sense of Bohr \cite{bohr}.
Let $d$ be the maximal number of linearly independent $\tau_1,...,\tau_m$ over $\mathbb{Q}.$ $A$ is periodic if and only if $d = 1.$
Furthermore, there exists a polynomial $P(z_1,...,z_d)$ and real numbers $r_1,...,r_d$ such that
\begin{equation}
\label{AP}
    A(y) = P(e^{2\pi i r_1y},...,e^{2\pi i r_dy}), \ \ y \in \mathbb{R}.
\end{equation}
In remainder of this section we derive the following two
theorems.
\begin{theo}
\label{theorem1}
$\ln |A|$ is almost periodic, in the $($mean squared$\,)$ sense of Besicovitch \cite{besicovitch}, and
\begin{equation}
\label{mean}
\lim_{L \rightarrow \infty} \frac{1}{2L} \int_{-L}^{L} \ln |A(y)| = \ln M(P).
\end{equation}
\end{theo}
\begin{defi}
\label{MA}
The Mahler measure of $A$ is $M(A) = M(P).$
\end{defi}
\begin{theo}
\label{theorem2}
    Define $\rho(f) = - \ln M(A)/\ln \lambda.$ Then
\begin{equation}
\label{Asympt1}
    \lim_{L \rightarrow \infty} \frac{1}{2L \ln L} \int_{-L}^{L} \ln |{\widehat f}\, (y)|\, dy = - \rho(f).
\end{equation}
\end{theo}
\begin{lem}
\label{titchmarsh}
$(${\bf Titchmarsh}$)$
If $h$ is a compactly supported distribution and $[a,b]$ is
the smallest interval containing the support of $h$ then
\begin{equation}
\label{exptype}
    H(z) = \int_{a}^{b} h(t) e^{2\pi zt} dt
\end{equation}
is an entire function of exponential type. Moreover the number $n(r)$
of zeros of $H(z)$ in the disk $|z| < r$ satisfies
\begin{equation}
\label{titchmarsheqn}
    \lim_{r \rightarrow \infty} \frac{n(r)}{r} = 2(b-a).
\end{equation}
\end{lem}
Titchmarsh's proof (Theorem IV, \cite{titchmarsh}) assumed that $h$ was a Lebesgue integrable function but a simple regularization argument implies that if $h$ is a compactly supported distribution then
Equation \ref{titchmarsheqn} holds. Therefore there exists $\beta > 0$ such that
\begin{equation}
\label{nzeros}
    \hbox{card } \{ \, y \, : \, y \in [-L,L] \hbox{ and } H(y) = 0 \, \}  \leq \beta\, L, \ \ L \geq 1
\end{equation}
\begin{lem}
\label{uvbound}
If $I = [a,b],$ $u > 0,$ $v > 0,$ $K : I \rightarrow \mathbb{C}$ is differentiable,
$|K| \leq v,$ the derivative $K^{\prime}$ satisfies
$|K^{\prime}| \geq u,$ and $K^{\prime}(I)$ is contained
in a single quadrant of $\mathbb{C},$ then
\begin{equation}
\label{muI}
    b - a \leq \frac{2{\sqrt 2}v}{u}.
\end{equation}
\end{lem}
{\bf Proof} We present the proof that we gave in (\cite{lawton3}, Lemma 1).
Since $u \leq |K^{\prime}|$ the triangle inequality
$|K^{\prime}| \leq |\Re\, K^{\prime}| + |\Im\, K^{\prime}|$
gives
$$u(b-a) \leq \int_{a}^{b} |K^{\prime}(y)| \, dy \leq \int_{a}^{b}
\left( \, |\Re \, K^{\prime}(y)| +  |\Im \, K^{\prime}(y)| \, \right) \, dy.$$
Since $K^{\prime}(I)$ is contained in a single quadrant of $\mathbb{C}$ there exist
$c \in \{1,-1\}$ and $d \in \{1,-1\}$ such that $|\Re \, K^{\prime}(y)| = c\, \Re \, K^{\prime}(y)$
and $|\Im \, K^{\prime}(y)| = d \, \Im \, K^{\prime}(y)$ for all $y \in I.$
Therefore
$$\int_{a}^{b}
\left( \, |\Re \, K^{\prime}(y)| +  |\Im \, K^{\prime}(y)| \, \right) \, dy
= (c\, \Re K(b) + d\, \Im K(b)) - (c\, \Re K(a) + d\, \Im  K(a)).$$
Since $|c| = |d| = 1$ and $|K(y)| \leq v$ the quantity on right
is bounded above by $2{\sqrt 2} \, v.$ Combining these three inequalities completes the proof.
\begin{lem}
\label{boyd}
Let $\mu$ denote Lebesque measure. There exists a sequence $c_2,c_3,...> 0$ such that for all monic polynomials $P(z)$ with $k \geq 2$ nonzero coefficients
\begin{equation}
\label{pb}
    \mu(\{\, y \in [0,1] \, : \, |P(e^{2\pi i \, y})| \leq v \, \}) \leq c_k\, v^{1/(k-1)}, \ \ v > 0.
\end{equation}
\end{lem}
{\bf Proof} We proved this in (Theorem 1, \cite{lawton3}) using induction on $k,$ lemma \ref{uvbound}, and the fact that a polynomial of degree $q$ can have at most $q$ distinct roots.
\begin{remark} We used Lemma \ref{boyd} and the method we developed in \cite{lawton2} and Boyd
developed in \cite{boyd} to prove a conjecture that he formulated in \cite{boyd}. The
conjecture expresses the Mahler measure of a multidimensional polynomial as a limit of
Mahler measures of univariate polynomials and provides an alternative proof of the
special case of Lehmer's conjecture proved in \cite{dobrowolski}. The conjecture is stated
in http:$//$en.wikipedia.org$/$wiki$/$Mahler$\textunderscore$measure
\end{remark}
\begin{lem}
\label{tb}
There exists a sequence $c_2,c_3,...> 0$ such that
\begin{equation}
\label{tb1}
    \mu(\{\, y \in [-L,L] \, : \, |A(y)| \leq v \, \})
    \leq L\, c_m \, v^{1/(m-1)}, \ \ L \geq 1, \ v > 0.
\end{equation}
\end{lem}
{\bf Proof} We observe that Lemma \ref{tb} concerns (not necessarily periodic) trigonometric polynomials. The conclusion is obvious for $m = 2$ so we assume that $m \geq 3$ and proceed by induction on $m$ using a technique similar to the one we used to prove Lemma \ref{boyd}. Since $|e^{-2\pi i y}A(y)| = |A(y)|$ we may assume without loss of generality that $0 = \tau_1 < \cdots \tau_m.$ Then $A^{\prime}(y)$ is a trigonometric polynomial with $m-1$ terms
so by induction there exists $c_{m-1} > 0$ such that
\begin{equation}
\label{tb2}
    \mu(\{\, y \in [-L,L] \, : \, |A^{\prime}(y)| \leq u \, \})
    \leq L\, c_{m-1} \, u^{1/(m-2)}, \ \ L \geq 1, \, u > 0.
\end{equation}
We observe that $A$ and $A^{\prime},$ their real and imaginary parts, and their squared moduli are restrictions of functions of exponential type so they satisfy the hypothesis in Lemma \ref{titchmarsh}. Therefore Equation \ref{nzeros} implies that there exists $\beta > 0$ such that for all $u >0, v > 0$ and $L \geq 1$ the subset of $[-L,L]$ where $|A^{\prime}(y)| \geq u$ and $|A(y)| \leq v$ can be expressed as the union of not more than $\beta\, L$ closed intervals $I$ such that $A^{\prime}(I)$ is contained in a single quadrant of $\mathbb{C}.$ Therefore Lemma \ref{uvbound} implies that there exists $\beta > 0$ such that
\begin{equation}
\label{tb3}
    \mu(\{\, y \in [-L,L] \, : \, |A(y)| \leq v \hbox{ and } |A^{\prime}(y)| \geq u \, \})
    \leq \frac{L\, \beta \, v}{u}, \ \ L \geq 1, \, u > 0, \, v > 0.
\end{equation}
Combining Equations \ref{tb2} and \ref{tb3} gives for every $L \geq 1,\, v > 0$
\begin{equation}
\label{tb4}
    \mu(\{\, y \in [-L,L] \, : \, |A(y)| \leq v  \, \}) \leq L \min_{u > 0}
        \left[\, c_{m-1} \, u^{1/(m-2)} + \frac{\beta \, v}{u} \, \right] = Lc_m v^{1/(m-1)}
\end{equation}
where
$c_{m} = (1+\beta(m-2))\left[\, \frac{c_{m-1}}{\beta(m-2)}\, \right]^{\frac{m-2}{m-1}}.$
\\ \\
For $n \in \mathbb{N},$ $\mathbb{R}^n$ is $n$-dimensional Euclidean vector space with norm $|| \, ||,$ $\mathbb{Z}^n$ is a rank $n$ lattice subgroup of $\mathbb{R}^n,$ $\mathbb{T}^n =  \mathbb{R}^n/\mathbb{Z}^n$ is the $n$-dimensional torus group equipped with Haar measure, and $p_n : \mathbb{R}^n \rightarrow \mathbb{T}^n$ is the canonical epimorphism. Vectors are represented as column vectors and $T$ denotes transpose. We define $|| \ ||_{\, \mathbb{T}^n}: \mathbb{T}^n \rightarrow \mathbb{R}$ by
\begin{equation}
\label{norm}
||\,g\,||_{\, \mathbb{T}^n} = \min \{ \, ||\,x\,|| : x \in \mathbb{R}^n \hbox{ and } p_n(x)=g \, \}, \ \ g \in \mathbb{T}^n.
\end{equation}
For $d \in \mathbb{N}$ and $r = [r_1,...,r_d]^T \in \mathbb{R}^d$ we define the homomorphism
$\Psi_r \, : \, \mathbb{R} \rightarrow \mathbb{T}^d$
by
\begin{equation}
\label{psi}
    \Psi_r(y) = p_d(y r), \ \ y \in \mathbb{R}
\end{equation}
where $p_d : \mathbb{R}^d \rightarrow \mathbb{T}^d$ is the canonical epimorphism.
\begin{lem}
\label{BSW}
$(${\bf Bohl--Sierpinski--Weyl}$)$ The image of $\Psi_r$ is a dense subgroup of $\mathbb{T}^d$ if and only if the components $r_1,...,r_d$ of $r$ are linearly independent
over $\mathbb{Q}.$ If they are independent then for every continuous function $S : \mathbb{T}^d \rightarrow \mathbb{C}$
\begin{equation}
\label{unifdist}
    \lim_{L \rightarrow \infty} \frac{1}{2L} \, \int_{-L}^{L} S \circ \Psi_r\, (y)dy = \int_{\mathbb{T}^d} S(g)\, dg
\end{equation}
where $S \circ \Psi_r : \mathbb{R} \rightarrow \mathbb{C}$ is the composition of $S$ with $\Psi_r$
and $dg$ is Haar measure on $\mathbb{T}^d.$
\end{lem}
{\bf Proof} Arnold (\cite{arnold}, p. 285--289) discuss the significance and gives Weyl's proof of
this classic result called the {\it Theorem on Averages} and asserts it ``may be found implicitly
in the work of Laplace, Lagrange, and Gauss on celestial mechanics'' and ``A rigorous proof was
given only in 1909 by P. Bohl, V. Sierpinski, and H. Weyl in connection with a problem of
Lagrange on the mean motion of the earth's perihelion."
\\ \\
For every nonzero polynomial $P(z_1,...\, ,z_d)$
we define the trigonometric polynomial \\
$P_{\hbox{trig}} : \mathbb{T}^d \rightarrow \mathbb{C}$
by
\begin{equation}
\label{Ptrig}
    P_{trig}(p_d([t_1,...\, ,t_d]^T)) = P(e^{2\pi i\, t_1},...\, ,e^{2\pi i\, t_d}).
\end{equation}
and observe through a simple induction based computation that
    $\ln |P_{trig}| \in L^2(\mathbb{T}^d).$
\\
Henceforth we assume that $r = [r_1,...\, ,r_d]^T \in \mathbb{R}^d$ and
$P(z_1,...\, ,z_d)$ are chosen to satisfy Equation \ref{AP}. Therefore
$r_1,...\, ,r_d$ are linearly independent over $\mathbb{Q}$ and moreover
\begin{equation}
\label{AP1}
    A(y) = P_{trig} \circ \Psi_r (y), \ \ y \in \mathbb{R}.
\end{equation}
{\bf Proof of Theorem \ref{theorem1}}
For all $v > 0$ define $S_v : \mathbb{T}^n \rightarrow \mathbb{R}$ by
$$
    S_v(g) = \ln \left( \max \{v,|P_{trig}(g)|\} \right).
$$
Lebesgue's dominated convergence theorem implies that
\begin{equation}
\label{LDC}
    \lim_{v \rightarrow 0} \int_{\mathbb{T}^d} |\, S_v(g)-\ln |P_{trig}(g)| \, |^2 dg = 0.
\end{equation}
Lemma \ref{tb} implies that
\begin{equation}
\label{Av}
    \lim_{v \rightarrow 0} \lim_{L \rightarrow \infty} \frac{1}{2L}
    \int_{-L}^{L} |\, \ln |A(y)| - S_v \circ \Psi_r(y)) \, |^2 dy = 0
\end{equation}
and hence $\ln |A(y)|$ is almost periodic in the sense of Besicovitch.
We combine Equations \ref{LDC} and \ref{Av} with Lemma \ref{BSW} to complete the proof
by computing
$$
\begin{array}{ccccc}
    \lim_{L \rightarrow \infty} \frac{1}{2L} \int_{-L}^{L} \ln |A(y)| dy & = &
    \lim_{v \rightarrow 0} \lim_{L \rightarrow \infty} \frac{1}{2L} \int_{-L}^{L} S_v \circ \Psi_r(y)) dy & = & \\ \\
    \lim_{v \rightarrow 0} \int_{\mathbb{T}^d} S_v(g) dg & = &
    \int_{\mathbb{T}^d} \ln |P_{trig}(g)| dg & = & M(P).
\end{array}
$$
\begin{remark} Since the zero set of $P_{trig}$ is a real analytic set, an alternative proof
based on Lojasiewicz's structure theorem \cite{kranzparks}, \cite{lojasiewicz} for real analytic sets may be possible using methods that we developed in \cite{lawton5} to prove the Lagarias-Wang Conjecture.
\end{remark}
{\bf Proof of Theorem \ref{theorem2}} Let
$L = b|\lambda|^k$ with $k \in \mathbb{N}$ and $b \in [\, |\lambda^{-1}|,1\, )$ then use the functional Equation \ref{functional} and Theorem \ref{theorem1} to compute
$$
\begin{array}{ccc}
\lim_{L \rightarrow \infty} \frac{1}{2L \ln L} \int_{-L}^{L} \ln |{\widehat f}\, (y)|\, dy
& = & \\ \\
\lim_{k \rightarrow \infty} \frac{1}{2 kb\, |\lambda|^k \ln |\lambda|}
\sum_{j = 1}^{k-1} \int_{b|\lambda|^j}^{b|\lambda|^{j+1}} \ln \left(\, |{\widehat f}\, (y){\widehat f}\, (-y)|\, \right)\, dy
& = & \\ \\
\lim_{k \rightarrow \infty} \frac{1}{2 k\, |\lambda|^k \ln |\lambda|}
\sum_{j = 1}^{k-1} |\lambda|^j \int_{1}^{|\lambda|}
\ln \left(\, |{\widehat f}\, (ub|\lambda|^j){\widehat f}\, (-ub|\lambda|^j)|\,\right)\,  du
& = & \\ \\
\lim_{k \rightarrow \infty} \frac{1}{2 k\, |\lambda|^k \ln |\lambda|}
\sum_{j = 1}^{k-1} |\lambda|^j \int_{1}^{|\lambda|}
\ln \left(\, |\, {\widehat f}\, (ub\lambda^j)\, {\widehat f}\, (-ub\lambda^j)\, |\,\right)\, du
& = & \\ \\
\lim_{k \rightarrow \infty} \frac{1}{2 k\, |\lambda|^k \ln |\lambda|}
\sum_{j = 1}^{k-1} |\lambda|^j \int_{1}^{|\lambda|} \left(\ln |{\widehat f}\, (bu){\widehat f}\, (-bu)| + \sum_{i=1}^j \ln |A(ub\lambda^j)A(-ub\lambda^j)|\right)du
& = & \\ \\
\lim_{k \rightarrow \infty} \frac{1}{2 k\, |\lambda|^k \ln |\lambda|}
\sum_{j = 1}^{k-1} |\lambda|^j \sum_{i=1}^j \int_{1}^{|\lambda|} \left(
\ln |A(ub\lambda^j)| + \ln |A(-ub\lambda^j)|\right)\, du
& = & \\ \\
\lim_{k \rightarrow \infty} \frac{(|\lambda| - 1)\, \ln M(A)}{k\, |\lambda|^k \, \ln |\lambda|}
\sum_{j = 1}^{k-1} j\, |\lambda|^j
& = & \\ \\
\lim_{k \rightarrow \infty} \frac{(|\lambda| - 1)\, \ln M(A)}{k\, |\lambda|^k \, \ln |\lambda|}
\left(\, \frac{|\lambda|-|\lambda|^{k+1}}{(|\lambda| - 1)^{-2}} +
\frac{k|\lambda|^{k-1}-1}{|\lambda|-1}\, \right)
& = & \\ \\
\frac{\ln M(A)}{\ln \lambda}.
\end{array}
$$
\section{Quasilattices}
\label{quasilattices}
If $\lambda$ is an algebraic integer with minimal polynomial
$\Lambda(z) = z^n + c_{n-1}z^{n-1} + \cdots + c_0, \ c_j \in \mathbb{Z}$
that has roots $\lambda = \lambda_1, \lambda_2, \dots, \lambda_n,$ then the
Frobenius companion matrix
\begin{equation}
\label{frobenius}
    C =
    \left[\begin{array}{ccccc}
        0 & 1 & 0 & \hdots & 0 \\
        0 & 0 & 1 & \ddots & \vdots \\
        \vdots & \ddots & \ddots & \ddots & 0 \\
        0 & 0 & \hdots & 0 & 1 \\
        -c_0 & -c_1 & \hdots & \hdots & -c_{n-1}
    \end{array}\right],
\end{equation}
Vandermonde matrix
\begin{equation}
\label{vandermonde}
    V =
    \left[\begin{array}{cccc}
     1                 &  1                & \hdots &  1 \\
     \lambda_1         & \lambda_2         & \hdots & \lambda_{n} \\
     \vdots            & \vdots            & \hdots & \vdots \\
     \lambda_{1}^{n-1} & \lambda_{2}^{n-1} & \hdots & \lambda_{n}^{n-1} \\
    \end{array} \right],
\end{equation}
and diagonal matrix $D = \hbox{diag}(\lambda_1,...\lambda_n)$ satisfy
\begin{equation}
\label{keyeq}
    C^k \, V = V \, D^k, \ \ k \in \mathbb{N}.
\end{equation}
A Pisot–-Vijayaraghavan (PV) number is an algebraic integer $\lambda$ whose conjugates $\lambda = \lambda_1,...,\lambda_n$ satisfy
$|\lambda_j| < 1$ whenever $j > 1.$
Then $|\lambda| > 1.$ They were discovered by Thue \cite{thue}, rediscovered by Hardy \cite{hardy} who with Vijayaraghavan studied their Diophantine approximation \cite{cassels} properties, and studied by Pisot \cite{pisot} in his dissertation. Examples are:
$$
\begin{array}{ccccc}
\lambda =\lambda_1  & \lambda_2 & \lambda_2 & \lambda_4 & \hbox{minimal polynomial} \\
1.6180 &-0.6180 &  &  &  z^2 - z - 1 \\
-1.6180 & 0.6180 & & & z^2 + z - 1 \\
3.4142 &  0.5858 &  &  &  z^2 - 4z + 2 \\
2.2470 & 0.5550 & -0.8019 & & z^3 - 2z^2 -z + 1 \\
1.3247 & -0.6624 + 0.5623i & -0.6624 - 0.5623i &  &  z^3-z-1\\
1.3803 &  -0.8192 & 0.2194 + 0.9145i  &  0.2194 - 0.9145i &  z^4-z^3 -1
\end{array}
$$
For $\lambda \in \mathbb{C}$ let $\mathbb{Z}[\lambda]$ be the $\mathbb{Z}$--module generated by powers of $\lambda.$
\begin{lem}
\label{pvnorm}
    If $\lambda$ is a PV number and $\alpha \in \mathbb{Z}[\lambda]$ then
    there exists a sequence $n_k \in \mathbb{Z}$ such that
    $\lambda^k \alpha - n_k \rightarrow 0$ exponentially fast.
\end{lem}
{\bf Proof} Define $n_k = \sum_{j=1}^{n} \lambda_j^k \alpha_j$
where $\lambda_1 = \lambda,...,\lambda_n$ and $\alpha_1 = \alpha,...,\alpha_n$ are the Galois conjugates of $\lambda$ and $\alpha.$ Each $n_k$ is a sum of the conjugates of $\lambda^k \alpha$ and hence is an integer. The assertion follows since $\sum_{j=2}^n \lambda_j^k \alpha_j \rightarrow 0$ exponentially fast.
\\ \\
Let $\lambda$ be a PV number and let $C, D, V,\lambda_1,...,\lambda_n$ be as above. We call a vector $\sigma \in \mathbb{R}^n$ admissible if
$\sigma_1 = 0,$ $\sigma_j > 0$ for $j \geq 2,$ and whenever $2 \leq j < k \leq n$
and $\lambda_j$ and $\lambda_k$ are nonreal complex conjugate pairs then $\sigma_j = \sigma_k.$
We denote the set of admissible vectors by $\mathbb{R}_{a}^{n}.$ It is an open cone and admits the
partial order $\sigma \leq \xi$ if and only if $\sigma_j \leq \xi_j,\ 1 \leq j \leq n.$
\begin{defi}
\label{quasilattice}
The quasilattice corresponding to $\sigma \in \mathbb{R}_{a}^n$ is
$$
  \mathfrak{L}(\sigma) =
  \{\, (V^T\ell)_1 \, : \, \ell \in \mathbb{Z}^n \hbox{ and } |(V^T\ell)_j| < \sigma_j \hbox{ for all } 2 \leq j \leq n \, \}.
$$
\end{defi}
\begin{remark} Quasilattices are a subset of cut-and-project sets studied by Meyer \cite{meyer}. de Bruijn \cite{debruijn}, \cite{auyangperk} showed that quasicrystals associated to Penrose tilings \cite{penrose1,penrose2,penrose3} are higher dimensional cut-and-project sets. Bombieri and Taylor \cite{bombieri1,bombieri2} discuss the diffraction and algebraic properties of quasicrystals. Hof \cite{hof} proved that every cut-and-project set has a pure point Fourier spectrum. The book by Senechal \cite{senechal} and references therein provide a basic introduction to quasicrystals.
\end{remark}
\begin{lem}
\label{QL1}
If $\sigma, \xi \in \mathbb{R}_{a}^n$ then
$$0 \in \mathfrak{L}(\sigma) \hbox{ and } \mathfrak{L}(\sigma) = - \mathfrak{L}(\sigma),$$
$$\mathfrak{L}(\sigma) \subseteq \mathfrak{L}(\xi) \hbox{ if and only if } \sigma \leq \xi,$$
$$\mathfrak{L}(\sigma) + \mathfrak{L}(\xi) = \mathfrak{L}(\sigma + \xi).$$
\end{lem}
{\bf Proof} Last two assertions follow from density property in Lemma \ref{BSW}, others are obvious.
\begin{lem}
\label{QL2}
Let $\sigma \in \mathbb{R}_{a}^n.$ If $\ell \in \mathfrak{L}(\sigma) \backslash \{0\}$ then
$|(V^T\ell)_1| \geq \prod_{j=2}^n \sigma_{j}^{-1}.$
The minimal distance between the points in $\mathfrak{L}(\sigma)$ is
$2^{1-n}\prod_{j=2}^n \sigma_{j}^{-1}.$
\end{lem}
{\bf Proof} If $\ell \in \mathbb{Z}^n \backslash \{0\}$ then
$\prod_{j=1}^n (V^T\ell)_j \in \mathbb{Z} \backslash \{0\}$ is nonzero and a symmetric polynomial in
$\mathbb{Z}[\lambda_1,...,\lambda_n]$ so is a nonzero integer, implying the first assertion. The second assertion follows since lemma \ref{QL1} implies
that the differences between points in $\mathfrak{L}(\sigma)$ are in $\mathfrak{L}(2\sigma).$
\\ \\
The following result, proved by Minkowsky in 1896 \cite{minkowski}, is the foundation of
the geometry of numbers (\cite{gruberlekkerkerker}, II.7.2, Theorem 1), (\cite{maninpanchishkin},
Theorem 4.7)
\begin{lem}
\label{minkowski}
$(${\bf Minkowski}$)$ If $X \subset \mathbb{R}^n$ is convex, $X = -X,$ and the volume of $X$ exceeds
$2^n$ then $X$ contains a point in $\mathbb{Z}^n\backslash \{0\}.$
\end{lem}
\begin{lem}
\label{QL3}
If $\sigma \in \mathbb{R}_{a}^n$ and $L > |\det V| \, \prod_{j=2}^n \sigma_{j}^{-1}$ then $\mathfrak{L}(\sigma)$ contains a nonzero point in the interval $(-L,L).$
\end{lem}
{\bf Proof} Let $X = (V^T)^{-1}\left[(-L,L) \times (-\sigma_1,\sigma_1) \cdots (-\sigma_n,\sigma_n)\right]^T.$ Then $X = -X$ and
$$\hbox{volume } X = 2^n \, |\det V|^{-1}\, L\, \prod_{j=2}^n \sigma_j > 2^n.$$
The result follows from Minkowski's Lemma \ref{minkowski} since
if $\ell \in \mathbb{Z}^n$ then
$(V^T\ell)_{1} \in (-L,L) \cap \mathfrak{L}(\sigma)$ if and only if
$
    \ell \in X.
$
\begin{lem}
\label{QL4}
If $\sigma \in \mathbb{R}_{a}^n$ then there exists $M(\sigma) > 0$ such that every open interval of length $M(\sigma)$ contains a point in $\mathfrak{L}(\sigma).$
\end{lem}
{\bf Proof} Clearly matrix $VV^T$ is invertible and it has integer entries since
\begin{equation}
\label{keyeq2}
    (VV^T)_{i,j} = \sum_{\ell = 1}^{n} \lambda_{\ell}^{i+j-2} = \hbox{trace } C^{i+j-2}, \ \ i,j \in \{1,...,n\}.
\end{equation}
Therefore $VV^T \mathbb{Z}^n$ is a rank $n$ subgroup of $\mathbb{Z}^n$ so the quotient group $G = \mathbb{R}^n/(VV^T\mathbb{Z}^n)$ is isomorphic to the $n-$dimensional torus group $\mathbb{T}^n.$ Let
$q_n : \mathbb{R}^n \rightarrow G$
and
$\varphi_n : G \rightarrow \mathbb{T}^n$
be the canonical epimorphisms and observe that
\begin{equation}
\label{equivariant}
    p_n\left((VV^T)^{-1}x\right) = \varphi_n \circ q_n(x), \ \ x \in \mathbb{R}^n.
\end{equation}
Since the entries of $v = [\, 1, \, \lambda, \, ... \, , \lambda^{n-1}\, ]^T$ are linearly independent over $\mathbb{Q}$ the density assertion in Lemma \ref{BSW} implies that
\begin{equation}
\label{density1}
    G = q_n\left( (V
    \left[\, \mathbb{R} \times (-\sigma_1,\sigma_1) \times (-\sigma_n,\sigma_n) \, \right]^T \right).
\end{equation}
Since $G$ is compact there exists $\kappa >0$ such that
\begin{equation}
\label{density2}
    G = q_n\left( (V
    \left[\, (-\kappa,\kappa) \times (-\sigma_1,\sigma_1) \times (-\sigma_n,\sigma_n) \, \right]^T \right).
\end{equation}
Therefore Equation \ref{equivariant} implies that
\begin{equation}
\label{density3}
    \mathbb{T}^n = \varphi_n(G) = p_n\left( (V^T)^{-1}
    \left[\, (-\kappa,\kappa) \times (-\sigma_1,\sigma_1) \times (-\sigma_n,\sigma_n) \, \right]^T \right)
\end{equation}
and the result follows by choosing $M(\sigma) = 2\kappa.$
\begin{remark} Lemma \ref{QL2} shows that quasilattices are uniformly
discrete, Lemma \ref{QL4} shows that they are relatively dense. Sets having
both of these properties are called Delone sets, named after Boris Nikolaevich Delone (aka Delaunay) \cite{delone}.
\end{remark}
An inflation symmetry of $\mathfrak{L}(\sigma)$ is a real number $\eta > 1$ such that $\eta \, \mathfrak{L}(\sigma) \subseteq \mathfrak{L}(\sigma).$
\begin{lem}
\label{QL5}
If $|c_0| = 1$ then
$\lambda \, \mathfrak{L}(\sigma) = \mathfrak{L}([\, 0,\, |\lambda_2|\, \sigma_2,...,\, |\lambda_n|\, \sigma_n\, ]^T).$ Therefore $\lambda$ is an inflation symmetry of $\mathfrak{L}(\sigma).$
\end{lem}
{\bf Proof}
Equation \ref{keyeq} gives $\lambda_j(V^T\ell)_j = (V^TC^T\ell)_j, \ j = 1,...,n$
and hence
$$
\lambda \mathfrak{L}(\sigma) =
    \{ \, (V^TC^T\ell)_1 \, : \, \ell \in \mathbb{Z}^n \hbox{ and }
        |(V^TC^T\ell)_j| < |\lambda_j|\, \sigma_j \hbox{ for } j = 2,...,n\, \}.
$$
Since $\det C = (-1)^n c_0 = \pm 1,$ $C^T \mathbb{Z}^n = \mathbb{Z}^n$ and hence
$\lambda \mathfrak{L}(\sigma) = $
$$
        \{ \, (V^T \, \ell)_1 \, : \, \ell \in \mathbb{Z}^n \hbox{ and }
        |(V^T\, \ell)_j| < |\lambda_j|\, \sigma_j \hbox{ for } j = 2,...,n\, \} =
        \mathfrak{L}([\, 0,\, |\lambda_2|\, \sigma_2,...\, ,|\lambda_n|\, \sigma_n\, ]^T).
$$
The second assertion then follows from Lemma \ref{QL1} since $|\lambda_j| < 1$ for $j = 2,...,n.$
\\ \\
Points $x < y$ in $\mathfrak{L}(\sigma)$ are consecutive if $x \neq y$ and $(x,y) \cap \mathfrak{L}(\sigma) = \phi.$ Define
$$\mathfrak{D}(\sigma) = \{\, y-x: x < y \hbox{ and } x \hbox{ and } y \hbox{ are consecutive points in } \mathfrak{L}(\sigma)\,  \}.$$
Lemma \ref{QL2} implies that $\mathfrak{D}(\sigma)$ is bounded below by a positive number and Lemma \ref{QL4} implies that $\mathfrak{D}(\sigma)$ is bounded above. The following result is stronger.
\begin{lem}
\label{QL6}
$\mathfrak{D}(\sigma)$ is finite.
\end{lem}
{\bf Proof} Let $x < y$ be consecutive points in $\mathfrak{L}(\sigma).$ Lemma \ref{QL1} implies that $y - x \in \mathfrak{L}(2\sigma)$ so Definition \ref{quasilattice} implies that there exists $\ell \in \mathbb{Z}^n$ such that
$y - x = (V^T\ell)_1$ and $|(V^T\ell)_j| < 2\sigma_j$ for $j > 1.$ Lemma \ref{QL4} implies that $|(V^T\ell)_1| = |y - x| < M(2\sigma).$ Therefore $V^T\ell$ and hence $\ell$ must belong to a bounded subset of $\mathbb{R}^n$ and therefore $\ell$ must belong to a finite subset of $\mathbb{Z}^n$ so $\mathfrak{D}(\sigma)$ is finite.
\begin{lem}
\label{QL7}
There exists a finite subset $F \subset \mathbb{R}$ such that
$\mathfrak{L}(\sigma) + \mathfrak{L}(\sigma) \subseteq \mathfrak{L}(\sigma) + F.$
\end{lem}
{\bf Proof} Let $x, y \in \mathfrak{L}(\sigma)$ and $z = x+y.$ Lemma \ref{QL4} implies there exists $w \in \mathfrak{L}(\sigma)$ such that
$|z-w| < M(\sigma).$ Lemma \ref{QL1} implies that $z$ and $w$ are in the Delone set $\mathfrak{L}(2\sigma)$ so there exists an integer $m$ bounded by a function of $\sigma$ and sequence of consecutive points $z_0 = w, z_1, ... , z_m = z$ in
$\mathfrak{L}(2\sigma).$ Lemma \ref{QL6} implies that $\mathfrak{D}(2\sigma)$ is finite and hence the set of possible values for
$z - w = (z_m-z_{m-1}) + \cdots + (z_1-z_0)$ lies in the finite $m-$fold sum set
$F = \mathfrak{D}(2\sigma) + \cdots + \mathfrak{D}(2\sigma).$
\\ \\
A Salem number is an algebraic integer $\lambda$ whose conjugates $\lambda = \lambda_1,...,\lambda_n$ satisfy
$|\lambda_j| \leq 1$ whenever $j > 1$ and $|\lambda| > 1.$
Like PV numbers, Salem numbers appear in Diophantine approximation and harmonic analysis \cite{meyer,salem1,salem2}. The smallest $1.17628$ is the Mahler measure of Lehmer's polynomial
$z^{10}+z^9-z^7-z^6-z^5-z^4-z^3+z+1$
and is conjectured to be the smallest known Mahler measure $> 1$ of any polynomial.
\begin{remark} Lemma \ref{QL7} shows that $\mathfrak{L}(\sigma)$ is a Meyer set \cite{meyer2} and Lagarias \cite{lagarias} showed that every inflation symmetry of a Meyer set is a Salem number.
\end{remark}
Let $\mathfrak{D}(\sigma) = \{c_1,...,c_m\}$ and $I_j = [0,c_j), j = 1,...,m$ where
without loss of generality we assume that $c_1 \in \mathfrak{L}(\sigma),$ and define $\mathfrak{I} = \{I_j:j=1,...,m\}.$ Then for every $j = 1,...,m$ there exists $T_j \subset \mathfrak{L}(\sigma)$ and a function $\Upsilon : T_j \rightarrow \mathfrak{I}$ such that $0 \in T_1$ and
$\Upsilon(0) = I_1$ and
\begin{equation}
\label{substitution}
    \lambda \, I_j = \bigcup_{\tau \in T_j} (\Upsilon(\tau) + \tau).
\end{equation}
It follows that for every $k \in \mathbb{N},$ the interval $[0,\lambda^kc_1)$ can be expressed as a disjoint union of translates of intervals in the set $\mathfrak{I}$ by elements in $\mathfrak{L}(\sigma)$ and that every positive element in $\mathfrak{L}(\sigma)$ arises as a translate for some $k.$ The negative elements arise similarly. This gives a way of constructing the quasilattice by a substitution dynamical system. Schlottmann \cite{schlottmann} and Sirvent and Wang \cite{sirvent} discuss the relationship between quasicrystals, tilings and substitution dynamical systems.
\section{Multiresolution Analyses}
\label{multiresolution}
We refer to the reader to the textbook by Resnikoff and Wells \cite{resnikoff} for background on multiresolution analyses and its relation to refinement equations (called scaling equations) and to multiwavelets. In this paper we only address the scaling functions and do not attempt to construct wavelets from them.
\begin{theo}
\label{multires}
If $|c_0| = 1$ and $\sigma \in \mathbb{R}_{a}^n$ then
$\xi = [\, 0,\, (1-|\lambda_2|)\, \sigma_2,...\, ,(1-|\lambda_n|)\, \sigma_n\, ]^T \in \mathbb{R}_{a}^n$
and if $f$ is a refinable distribution satisfying Equation \ref{refinable}
and $\tau_j \in \mathfrak{L}(\xi)$ then the spaces
\begin{equation}
\label{W}
    W_k = \hbox{span } \{\, f(\lambda^k\, x - \tau)\, : \, \tau \in \mathbb{L}(\sigma)\, \},
    \ \ k \in \mathbb{Z},
\end{equation}
satisfy $W_k \subset W_{k+1}.$
\end{theo}
{\bf Proof} If $\tau \in \mathbb{L}(\sigma)$ then Equation \ref{refinable} implies
that
$$f(\lambda^k \, x - \tau) = \sum_{j=1}^m a_j f(\lambda^{k+1} \, x - \lambda \tau - \tau_j).$$
Since $|c_0| = 1,$ Lemma \ref{QL5} implies that
$$
\lambda \, \tau \in \mathfrak{L}([\, 0,\, |\lambda_2|\, \sigma_2,...,\, |\lambda_n|\, \sigma_n\, ]^T),
$$
and hence Lemma \ref{QL1} implies that
$$
\lambda\, \tau + \tau_j \in
\mathfrak{L}([\, 0,\, |\lambda_2|\,\sigma+\xi_2,...\, ,|\lambda_n|\, \sigma_n+\xi_n\, ]^T) = \mathfrak{L}(\sigma).
$$
Theorem \ref{multires} provides a multiresolution analysis for spaces of distributions. We have been unable to construct the refinable distribution $f$ that is an integrable function. In the next section we conjecture that such a construction is impossible.
\\ \\
We now construct a vector valued integrable function $f$ that is refinable meaning that it satisfies Equation \ref{refinable} when the coefficients $a_j$ are suitable matrices. $f$ is called a vector valued scaling function in the multiwavelet discussion in (\cite{resnikoff}, p. 137--138). We choose the entries of $f$ to be the characteristic functions $\chi_{I_j}, j = 1,...,m$ of the intervals in $\mathfrak{I}.$ Then Equation \ref{substitution} gives
\begin{equation}
\label{vector}
    \chi_{I_j}(x) = \sum_{\tau \in T_j} \chi_{\Upsilon(\tau)}(\lambda x - \tau), \ \ j = 1,...,m
\end{equation}
which determines the matrices $a_j.$ Suitable translates of entries of $f$ can be used to define the multiresolution spaces $W_k$ so that $W_k$ equals the set of functions that are constant on intervals of the form $[x,y)$ where $x$ and $y$ are consecutive points in $\lambda^{-k}\mathfrak{L}(\sigma).$ Refinable vectors with smoother entries are constructed from convolutions of entries in $f.$
\section{Conjectures}
\label{conjectures}
The Hausdorff dimension of a Borel measure is defined
in (\cite{peresschlagsolomyak}, p. 6) as
\begin{equation}
\label{HD}
    \hbox{dim}(\nu) = \inf \{ \, \hbox{dim}(B) \, : \, B \hbox{ is a Borel subset of } \mathbb{R} \hbox{ and }
    \nu (\mathbb{R} \backslash B) = 0 \ \}
\end{equation}
where $\hbox{dim}(B)$ is the Hausdorff dimension of $B.$
\begin{conj}
\label{hausdim}
If a refinable distribution $f$ is a Borel measure then $\rho(f) = \hbox{dim}(f).$
\end{conj}
If conjecture \ref{hausdim} is validated we propose to call $\rho(f)$ defined in Theorem \ref{theorem2} the Hausdorff dimension of $f.$ The following observations support this conjecture.
\\ \\
The characteristic function $f$ of $[0,1]$ satisfies the refinement equation $f(x) = f(2x) + f(2x-1).$ Therefore
$A(y) = (1 + e^{2\pi i y})/2$ so $\rho(f) = 1$ which equals the Hausdorff dimension of $f.$ We observe that $|{\widehat f}\, (y)| = |\sin y / y|$ decays like $|y|^{-1}.$
\\ \\
The uniform measure $f$ on Cantor's ternary set satisfies the refinement equation $f(x) = \frac{3}{2}f(3x) + \frac{3}{2}f(3x-2).$ Therefore $A(y) = (1 + e^{6\pi i y})/2$ so $\rho(f) = \ln 2/\ln 3 \approx 0.6309$ which equals the Hausdorff dimension of $f.$ This suggests that ${\widehat f}(y)$ decays like $|y|^{-\ln 2/\ln 3}.$
\\ \\
The measure $f$ defined by a Bernoulli convolution and PV number $\lambda$ satisfies the refinement equation
$f(x) = \frac{|\lambda|}{2}f(\lambda x) + \frac{|\lambda|}{2}f(\lambda x - 1).$
Therefore $A(y) = (1 + e^{2\pi i \lambda y})/2$ so $\rho(f) = \ln 2/\ln |\lambda|$ which equals the Hausdorff dimension of $f$ computed by Peres, Schlag and Solomnyak \cite{peresschlagsolomyak}. This suggests that
${\widehat f}(y)$ decays like $|y|^{-\ln 2/\ln |\lambda|}.$
\\ \\
Let $f$ be a nonzero distribution that satisfies the refinement equation $f(2x) = 2f(2x) + 2f(2x-1) - 2f(2x-3).$
Therefore $A(y) = 1+e^{2\pi i y} - e^{4\pi i y}$ so $\rho(f) = -\ln ((1+\sqrt 5)/2)/\ln 2 \approx 0.6942.$ This suggests that $|{\widehat f}(y)|$ grows like $|y|^{-\rho(f)}$
\\ \\
If $\lambda$ is a PV number, $m = 2,$ $a_1 = a_2 = |\lambda|/2,$ $\tau_1 = 0,$ and $\tau_2 = 1$ then $f$ is a measure defined by a Bernoulli convolution and Erd\"{o}s \cite{erdos} used Lemma \ref{pvnorm} to prove that there exists $\alpha > 0$ and $\gamma \in \mathbb{C} \backslash \{0\}$ such that
\begin{equation}
\label{erdos}
    \lim_{k \rightarrow \infty} {\widehat f}\, (\alpha \lambda^k) = \gamma.
\end{equation}
This proves that $f$ is not integrable since otherwise the Riemann-Lebesgue lemma implies that the limit in Equation \ref{erdos} equals $0.$
\\ \\
Kahane \cite{kahane} extended Erd\"{o}s' result for $\lambda$ a Salem number and Dai, Feng and Wang \cite{daifengwang} extended Kahane's result for $m > 2$ and $\tau_j \in \mathbb{Z}.$
\begin{conj}
\label{bigconj}
If $f$ is a nonzero refinable distribution with respect to dilation by a PV number $\lambda$ and translations in $\mathbb{Z}[1,\lambda,...,\lambda^{n-1}]$ then there exists $\alpha \in \mathbb{Q}[\lambda]\, \backslash \, \{0\}$ such that
$$\lim_{k\rightarrow \infty} {\widehat f}\, (\lambda^k\, \alpha) = \gamma$$
where
$\gamma \in \mathbb{C}\, \backslash \, \{0\}.$ This implies, by the Riemann--Lebesgue lemma, that $f$ is not integrable
thus extending the results in \cite{daifengwang} to certain noninteger values of $\tau_j.$
\end{conj}
A proof of this conjecture eludes us but we present some observations that support it. Henceforth let $\lambda$ be a PV number of degree $n,$ let $C : \mathbb{T}^n \rightarrow \mathbb{T}^n$ be the endomorphism induced by its Frobenius matrix defined in Equation \ref{frobenius}, let $f$ be a refinable distribution that satisfies Equation \ref{refinable} with $\tau_j \in \mathbb{Z}[\lambda],$ let $p_n : \mathbb{R}^n \rightarrow \mathbb{T}^n$ be the canonical epimorphism, and let $P_{trig} : \mathbb{T}^d : \rightarrow \mathbb{C}$ be the trigonometric polynomial that satisfies
\begin{equation}
\label{tp1}
    A(y) = P_{trig}\left(p_n([y,\lambda y,...,\lambda^{n-1}y]^T)\right), \ \ y \in \mathbb{R}.
\end{equation}
Therefore Equation \ref{keyeq} implies that
\begin{equation}
\label{tp2}
    A(\lambda \, y) = P_{trig}\left(C(p_n([y,\lambda y,...,\lambda^{n-1}y]^T))\right), \ \ y \in \mathbb{R}.
\end{equation}
Define $V(P_{trig}) = \{x \in \mathbb{T}^n \, : \, P_{trig}(x) = 0\, \}.$
\begin{lem}
\label{lemmaconjsupp}
     If $f$ is integrable and $\alpha \in \mathbb{Q}[\lambda]\, \backslash \, \{0\}$ then there exists an integer $k$ such that $A(\lambda^k\, \alpha) = 0.$ For $k \geq 0$ this is equivalent to
     \begin{equation}
     \label{super1}
         C^k\left(p_n([\alpha,\lambda \alpha,...,\lambda^{n-1} \alpha]^T)\right) \in \, V(P_{trig}),
     \end{equation}
     and for $k < 0$ this is equivalent to
     \begin{equation}
     \label{super2}
         p_n([\alpha,\lambda \alpha,...,\lambda^{n-1} \alpha]^T) \in \, C^{-k} \left(V(P_{trig})\right).
     \end{equation}
     If $|\det C| = |c_0| = 1$ then
     $C : \mathbb{T}^n \rightarrow \mathbb{T}^n$
     is invertible and Equation \ref{super2} is equivalent to
     \begin{equation}
     \label{super3}
         C^k\left(p_n([\alpha,\lambda \alpha,...,\lambda^{n-1} \alpha]^T)\right) \in \, V(P_{trig}).
     \end{equation}
\end{lem}
{\bf Proof} Lemma \ref{pvnorm} implies that
$\lim_{j \rightarrow \infty} \prod_{m = j}^{\infty} A(\lambda^m \alpha) = 1.$ Since $f$ is integrable the Riemann--Lebegue lemma implies that ${\widehat f}(\lambda^j \, \alpha) \rightarrow 0$ so the first assertion follows from Equation \ref{product}. The equivalences follow easily from Equation \ref{tp2} and the fact that $|\det C| = 1$ implies that the entries of $C^{-1}$ are integers.
\begin{remark}
Lemma \ref{lemmaconjsupp} shows that if $f$ is integrable then the real algebraic set $V(P_{trig})$ contains a rich set of points which suggests that it is not a proper subset of $\mathbb{T}^n$ and hence that Conjecture \ref{bigconj} holds.
\end{remark}
The following discussion shows that $V(P_{trig})$ may contain an even richer set of points than Lemma \ref{lemmaconjsupp} implies.
Let $A(\lambda)$ denote the $\mathbb{Q}-$subspace of $\mathbb{C}^n$ consisting of vectors $[\alpha_1,...,\alpha_n]^T$ whose components are conjugates and satisfy  $\alpha_j \in \mathbb{Q}(\lambda_j), \, j = 1,...,n.$ Then the map
$\Gamma_V \, : \, \mathbb{Q}^n \rightarrow A(\lambda)$ defined by
$\Gamma_V(q) = V^T(VV^T)^{-1}q, \ \ q \in \mathbb{Q}^n$
is a $\mathbb{Q}-$linear bijection of $\mathbb{Q}^n$ onto $A(\lambda)$ and
\begin{equation}
\label{keyeq3}
    q = V \Gamma_V q.
\end{equation}
Now assume that $\lambda$ is a PV number and let $v$ be the first column vector of $V$ and $\alpha_1$ be the first entry of $\Gamma_V q.$ Then Equations \ref{keyeq} and \ref{keyeq2} give
\begin{equation}
\label{conv1}
    ||\, C^k q - \alpha_1 \, \lambda^k \, v|| \rightarrow 0
\end{equation}
where the convergence is exponentially fast. Then Equation \ref{conv1}
implies that
\begin{equation}
\label{conv2}
    ||\, p_n(C^k q - \alpha \, \lambda^k \, v)||_{\mathbb{T}^n} = ||\, C^k p_n(q) - p_n(\alpha \, \lambda^k \, v)||_{\mathbb{T}^n} \rightarrow 0.
\end{equation}
Since $p_n(\mathbb{Q}^n)$ is the set of preperiodic points of the endomorphism $C : \mathbb{T}^n \rightarrow \mathbb{T}^n,$ if $\alpha \in \mathbb{Q}(\lambda)$ then the sequence $p_n(\alpha \, \lambda^k \, v) \in \mathbb{T}^n$ converges exponentially fast to a periodic orbit of
$C.$ Conversely every periodic orbit is the limit of such a sequence.
\begin{lem}
\label{lemmalast}
If $\lambda$ is a PV number with degree $n,$ $\alpha \in \mathbb{Q}(\lambda),$ $\{\tau_1,...,\tau_m\} \subset \mathbb{Z}\left[\, 1,\lambda,...,\lambda^{n-1}\, \right],$ and ${\widehat f}\, (\alpha \, \lambda^k) \neq 0$ for all $k \in \mathbb{N}$ then the following limit exists
\begin{equation}
\label{Asympt2}
    \lim_{k \rightarrow \infty} \frac{\ln|{\widehat f}\, (\alpha \, \lambda^k)|}{k\, \ln |\lambda|}
\end{equation}
and equals the mean value of $\ln |P_{\hbox{trig}}|$ over a finite orbit of $C$ in the $n$-dimensional torus group.
\end{lem}
{\bf Proof} Follows from the same calculation used to prove Theorem \ref{theorem2}.
\begin{remark}
Lemmas \ref{lemmaconjsupp} and \ref{lemmalast} ensure that $V(P_{trig})$ contains a rich set of points whenever $f$ is integrable. We think that the techniques that we used to prove the Lagarias--Wang conjecture in \cite{lawton5}, suitably modified, can be applied to prove Conjecture \ref{bigconj} by showing that there do not exist proper real algebraic subsets of $\mathbb{T}^n$ that contains these rich set of points. To this end we suggest studying the preperiodic points of $C$ using results from integral matrices discussed by Newman \cite{newman}, the open dynamical version of the Manin-Mumford Conjecture solved by Raynaud \cite{raynaud1}, \cite{raynaud2}, and connections between Mahler measure and torsion points developed by Le \cite{le}.
\end{remark}
\section{Multidimensional Extensions}
\label{multidimensional}
Constructing quasilattices in dimension $d > 1$ is simple. Choose any irreducible monic integral polynomial $\Lambda(x)$ of degree $n > 2$ and constant term $\pm 1$ having roots $\lambda_1,..,\lambda_d$ with modulus $> 1$ and whose remaining roots have modulus $< 1.$ The quasilattice $\mathfrak{L}(\sigma)$ consists of all vectors
$$[(V^T\ell)_1,...,(V^T\ell)_d]^T \in \mathbb{R}^d$$
where $\ell \in \mathbb{Z}^n$ are chosen to satisfy constraints of the form
$|(V^T\ell)_j| < \sigma_j$ for $j = d+1,..,n$ where $\sigma$ is an appropriate
admissible vector with $\sigma_1 = \cdots = \sigma_d = 0$ and the remaining entries positive. These multidimensional quasilattices have similar properties to the one dimensional quasilattices except that multiplication by $\lambda$ is replaced by multiplication by the diagonal matrix $\hbox{diag}[\lambda_1,...,\lambda_d].$
The construction of refinable distributions is similar but deriving their analytic properties is more difficult. It is easy to show the number of Voronoi cells
for these quasilattices is finite and and hence a finite number of dual Delaunay tiles, each a polytope whose vertices are in the quasilattice, can be constructed which tile $\mathbb{R}^d.$ However, although the quasilattice satisfies an inflation related to the Frobenius matrix for $\Lambda,$ in general the Delaunay tiling will not admit this inflation. One challenge is to determine possible symmetry conditions on $C$ that will allow this. We also suggest studying extensions obtained by replacing $\mathbb{R}^n$ by certain Lie groups $($stratified nilpotent groups with rational structure constants$)$ considered in \cite{lawton4}.

\Finishall
\end{document}